\documentclass[twocolumn]{autart}    

\usepackage{graphicx}          
\usepackage{amsmath}
\usepackage{amssymb}
\usepackage{enumitem}

\usepackage{savesym}
\savesymbol{AND}
\usepackage{algorithmic}

\usepackage{newfloat}
\DeclareFloatingEnvironment[name=Algorithm]{algofloat}
\newenvironment{algorithmn}[1]
{%
\begin{algofloat}[h]
\begin{minipage}{0.95\columnwidth}
\hrule width \hsize height 2pt \kern 1mm  \hrule width \hsize \kern 0.0em
\caption{#1}
\vspace{0.5em}
\hrule
}
{%
\hrule width \hsize \kern 1mm \hrule width \hsize height 2pt 
\end{minipage} 
\end{algofloat}
}

\usepackage{color}
\definecolor{gray90}{gray}{0.9}
\definecolor{gray80}{gray}{0.8}
\definecolor{gray70}{gray}{0.7}

\usepackage{cite}

\begin{document}

\begin{frontmatter}

\title{The Maximal Positively Invariant Set: Polynomial Setting}
\thanks[footnoteinfo]{Corresponding author S.~V.~Rakovic}

\author[Indep]{Sa\v{s}a.~V.~Rakovi{\'c}}\ead{sasa.v.rakovic@gmail.com},    
\author[ShTech]{Mario E.~Villanueva}\ead{meduardov@shanghaitech.edu.cn}               

\address[Indep]{Independent Researcher}  
\address[ShTech]{School of Information Science and Technology, 
                 ShanghaiTech University, Shanghai 201210, China.}  

\begin{keyword}                           
Set Invariance, Maximal Positively Invariant Set, Polynomial Dynamical Systems.
\end{keyword}                             

\begin{abstract}                          
This note considers the maximal positively invariant set for polynomial discrete time dynamics subject to constraints specified by a basic semialgebraic set. The note utilizes a relatively direct, but apparently overlooked, fact stating that the related preimage map preserves basic semialgebraic structure. In fact, this property propagates to underlying set--dynamics induced by the associated restricted preimage map in general and to its maximal trajectory in particular. The finite time convergence of the corresponding maximal trajectory to the maximal positively invariant set is verified under reasonably mild conditions. The analysis is complemented with a discussion of computational aspects and a prototype implementation based on existing toolboxes for polynomial optimization. 
\end{abstract}

\end{frontmatter}

\section{Introduction}
The theory of set invariance has been an important topic of an extensive research over the last five decades, due to its intimate relationships with fundamental notions and concepts within control and systems theory as well as theory of dynamical systems. For instance, set invariance analysis facilitates necessary and sufficient conditions for solvability of infinite horizon control under constraints and uncertainty synthesis problems. Furthermore, set invariance notions also lead to  necessary and sufficient conditions for maintaining deterministic or uncertain dynamics within prescribed constraints indefinitely long. Additional classical applications of set invariance include the design of stabilizing nominal/robust/stochastic optimal and model predictive controllers. The interested reader is referred 
to~\cite{bertsekas:1972,gilbert:tan:1991,aubin:1991,kolmanovsky:gilbert:1998,blanchini:miani:2008,artstein:rakovic:2008,
rawlings:mayne:2009,artstein:rakovic:2011} and numerous references therein for a related overview demonstrating that set invariance has earned an important place across a wide spectrum of classical fields including, \emph{inter alia}, backward and forward reachability analyses, deterministic/robust/stochastic optimal and model predictive control synthesis as well as stability analysis. Contemporary research encapsulates analysis and synthesis of smart autonomous systems as well as cyber--physical systems for which it is of paramount importance to guarantee operations with \emph{a--priori} guaranteed safety, resilience, security, reconfigurability  and fault tolerance. Such highly desirable \emph{modus operandi} can be ensured by adequate utilization of classical set invariance notions. Indeed, it should not be surprising that the maximal safe sets in contemporary jargon are equivalents of the maximal positively invariant sets in classical terminology. Furthermore, the design of resilient and secure systems demands a high degree of flexibility and robustness to uncertainty that can be offered by using robust control/positively invariant sets. Likewise, reconfigurability and fault tolerance is only possible from certain sets of states. A sample illustration of the use of positively invariant sets for fault tolerance can be found, for instance, in~\cite{olaru:dedona:seron:stoican:2011}. All in all, set invariance has a huge potential to significantly shape and affect analysis and synthesis of contemporary systems.

An important class of dynamical systems that has recently received a lot of attention is that induced by polynomial state transition maps. Even though polynomial dynamics are ubiquitous in science and engineering, their real power stems from the Weierstrass's approximation theorem~\cite{rudin:1964}, which in essence states that the set of polynomials is dense in the set of continuous functions (over compact domains). In turn, polynomial dynamics can be used to approximate as closely as desired more complex continuous nonlinear dynamics. Naturally, the same conclusions extend to constraints specified by basic semialgebraic sets, i.e. constrains defined as the solution sets of a finite number of polynomial inequalities (and equalities). Indeed, the class of basic semialgebraic sets is rather rich, and it encapsulates frequently encountered classes of polytopic and ellipsoidal sets. 

Some of the prime considerations in set invariance are the characterization and computation of the maximal positively invariant set~\cite{bertsekas:1972,gilbert:tan:1991,aubin:1991,kolmanovsky:gilbert:1998,blanchini:miani:2008}.  The topological properties and finite time computability of the maximal positively invariant set are particularly important questions for obvious reasons. In regards to computational aspects, conference contributions~\cite{liu:zhan:zhao:2011emsoft} and~\cite{li:liu:2016cacsd} consider, respectively, quantifier elimination methods, based on the Lie algerba, for testing whether a semialgebriac set is an 
invariant of a polynomial dynamical system as well as interval arithmetics approaches to computation of the maximal control/positively invariant sets for nonlinear switched dynamics. A recent article~\cite{korda:henrion:jones:2014} explores utilization of an infinite dimensional linear program (formulated in the space of occupational measures) for the computation of the maximal positively invariant set for polynomial dynamics subject to semialgebraic state and control constraints. This work also considers computation of approximations (that are not necessarily positively invariant) of the maximal positively invariant set by   constructing a hierarchy of linear matrix inequalities and dual sum--of--square problems.
In the computational context and with the polynomial setting in mind, a relevant problem---from the set invariance perspective---would be to determine whether the maximal positively invariant set for given polynomial dynamics and basic semialgebraic constraint set belongs to the class of basic semialgebraic sets. It is reasonably well known that, under certain regularity conditions, the maximal positively invariant set for given dynamics and constraint set can be computed as the topologically compatible limit of the maximal trajectory of the particular set--dynamics, whose state transition map is the preimage map of the original dynamics, restricted by the constraint set. In this note, we provide conditions ensuring the affirmative answer to the former question by exploiting a fundamental fact about basic semialgebraic sets that has apparently been overlooked within the context of maximal positively invariant sets. This fact essentially asserts that the class of basic semialgebraic sets is closed under the polynomial preimage and intersection operations. 

In this note, we provide a reasonably simple demonstration that basic semialgebraic sets are closed under the related restricted preimage maps. This fact is then used, in conjunction with mathematical induction, to show that the maximal trajectory of the associated set--dynamics is a sequence of basic semialgebraic sets. Finally, sufficient conditions are provided under which the associated maximal positively invariant set is itself a basic semialgebraic set.  A prototype algorithm for the computation of the maximal positively invariant set of given polynomial dynamics and basic semialgebraic constraint set is also outlined. This prototype algorithm operates on polynomials and can be implemented using existing tools for polynomial optimization~\cite{parrilo:2003}. 

\textbf{Paper Structure:}
Section~\ref{sec:02} describes the problem of our interest. Section~\ref{sec:03} establishes the key facts of our interest. Section~\ref{sec:04} discusses a prototype algorithm for the computation of maximal positively invariant sets in the polynomial/basic semialgebraic setting, while Section~\ref{sec:05} reports an academic and illustrative example. 
Section~\ref{sec:06} concludes the paper.

\textbf{Notation and Conventions:}
The sets of nonnegative integers and real numbers are denoted by $\mathbb{N}$ and $\mathbb{R}$. The set of $n$--dimensional real vectors is denoted by $\mathbb{R}^{n}$. The $j^\text{th}$ component of an $n$--dimensional 
vector $y$ is denoted by $y_{j}$. Whenever a variable is indexed, e.g. 
$y_{i}\in\mathbb{R}^{n}$, its $j^\text{th}$ component is written as $y_{(i,j)}$.
Moreover, for a vector $y\in\mathbb{R}^{n}$, its transpose is denoted by 
$y^{\intercal}$.
We write $g$ or $g(\cdot)$ for a function, and $g(y)$ for its value at a point 
$y$ in its domain. Let $g\ :\ \mathbb{R}^{n}\to\mathbb{R}^{n}$, $g^k$ denote the 
$k^\text{th}$ functional power of $g$, i.e. 
\begin{equation*}
g^{k}:=\underbrace{g\circ\cdots\circ g}_{k \text{ copies of } g}.
\end{equation*}
Consistently, $g^{k}(y)$ denotes the $k^\text{th}$ iterate of the function $g$ at point $y$ in its domain. For a $y\in\mathbb{R}^{n}$, the set of $n$--variate polynomials in $y$ with coefficients in $\mathbb{R}$ is denoted by $\mathbb{R}[y]$. For a $y\in\mathbb{R}^{n}$ and $g:\mathbb{R}^{n}\to\mathbb{R}^{m}$, we write  
$g\in\mathbb{R}^{m}[y]$ if and only if $g_{i}\in\mathbb{R}[y]$ for each $i\in\{1,\ldots,m\}$.

Herein we work with nonempty sets unless otherwise
stated. Also, when referring to maximality of sets and sets in a sequence of sets
we refer to maximality and term--wise maximality w.r.t. set inclusion. For clarity of presentation, we provide proofs of less obvious statements in the form of appendices.

\section{Problem Description}
\label{sec:02}
We consider discrete time dynamics given by
\begin{equation}
\label{eq:02.01}
x^{+} = f\left(x\right),
\end{equation}
where $x\in\mathbb{R}^{n}$ and $x^{+}\in\mathbb{R}^{n}$ denote, respectively,  
the current and successor states. The function $f\ :\ \mathbb{R}^{n}\to\mathbb{R}^{n}$ 
denotes the state transition map. Thus, at any current time $k\in\mathbb{N}$, 
the successor state is given by $x_{k+1} = f(x_{k})$. The state variables are 
subject to hard constraints (that apply to all time instances $k\in\mathbb{N}$)
\begin{equation}
\label{eq:02.02}
x \in \mathbb{X} \subset \mathbb{R}^{n}.
\end{equation}

The dynamics under consideration are polynomial.
\begin{assum}
\label{ass:02.01} 
The state transition map is a, known exactly, polynomial function, i.e. $f\in\mathbb{R}^{n}[x]$.
\end{assum}
The constraint set belongs to the class of basic semialgebraic sets.
\begin{defn}[Basic Semialgebraic Set]
A set $X$ in $\mathbb{R}^{n}$ is said to be a basic semialgebraic if it is the 
solution set of a finite number of polynomial inequalities.
\end{defn}
In fact, we work with compact, basic semialgebriac constraint sets. 
\begin{assum}
\label{ass:02.02} 
The constraint set $\mathbb{X}$ is compact and 
basic semialgebraic. It is given by
\begin{equation}
\label{eq:02.03}
\mathbb{X} := \left\{ x\in\mathbb{R}^{n}\ :\
\varphi_{0}( x )\geq 0 \right\}
\end{equation}
for some known exactly $\varphi_{0}\in\mathbb{R}^{m}[x]$.
\end{assum}

The notions of positively and maximal positively invariant sets are standard, but are recalled for the sake of completeness. 
\begin{defn}[Positively Invariant Set]
\label{def:02.01}
A set $\Omega$ in $\mathbb{R}^{n}$ is positively invariant for dynamics of~\eqref{eq:02.01} and constraint set of~\eqref{eq:02.02} if 
and only if $\Omega\subseteq\mathbb{X}$ and for all 
$x\in\Omega$ it holds that $f(x)\in\Omega$.
\end{defn}

\begin{defn}[Maximal Positively Invariant Set]
\label{def:02.02}
A set $\Omega_{\infty}$ in $\mathbb{R}^{n}$ is the maximal positively 
invariant set for dynamics of~\eqref{eq:02.01} and 
constraint set of~\eqref{eq:02.02} if and only if $\Omega_{\infty}$ is a positively invariant set for dynamics of~\eqref{eq:02.01} and  constraint set of~\eqref{eq:02.02} and it contains all positively invariant sets for dynamics of~\eqref{eq:02.01} and  constraint set of~\eqref{eq:02.02}.
\end{defn}

For typographical convenience, we alternate the terms ``a positively invariant
set" and ``the maximal positively invariant set" with the related complete expressions as specified in Definitions~\ref{def:02.01} and ~\ref{def:02.02}, respectively.
Since the state transition map $f$ and constraint set $\mathbb{X}$ are fixed 
throughout this note, no confusion should arise. 

Generally speaking, this note focuses on the topological structure of the maximal positively invariant set for dynamics of~\eqref{eq:02.01} and  constraint set of~\eqref{eq:02.02} as well as on related computational aspects. Speaking more precisely, this note demonstrates that in the computationally most relevant case---which occurs under reasonably mild and natural assumptions---the maximal positively invariant set is finitely determined, compact, and basic semialgebraic.

\section{Key Facts}
\label{sec:03}
The maximal positively invariant set for dynamics of~\eqref{eq:02.01} and 
constraint set of~\eqref{eq:02.02} is the topologically compatible limit of the maximal trajectory of
the set--dynamics specified by 
\begin{equation}
\label{eq:03.01}
{X}^{+} = {F}^{-1}_{\mathbb{X}}\left( {X} \right),
\end{equation}
where ${F}^{-1}_{\mathbb{X}}(\cdot)$ denotes the  preimage map of a set ${X}$ in $\mathbb{R}^{n}$ under $f$ restricted to $\mathbb{X}$, i.e.
\begin{align}
{F}^{-1}_{\mathbb{X}}({X}) &= \left\{ x\in\mathbb{R}^{n} \ : \
f(x)\in{X} \right\} \cap \mathbb{X} \nonumber\\
\label{eq:03.02}
&= \left\{ x\in\mathbb{X} \ : \ f(x)\in{X} \right\}.
\end{align}
More precisely, the related maximal positively invariant set  can be constructed~\cite{bertsekas:1972,gilbert:tan:1991,aubin:1991,kolmanovsky:gilbert:1998,blanchini:miani:2008} as the limit of the following set recursion
\begin{equation}
\label{eq:03.03}
\forall k\in\mathbb{N}, \quad {X}_{k+1} = {F}^{-1}_{\mathbb{X}}({X}_{k})\quad 
\text{with} \quad {X}_{0} = \mathbb{X}.
\end{equation}

\subsection{Structural Properties}
\label{sec:03.01}

We make use of a relatively direct, but highly relevant, fact that has apparently remained unexploited within set invariance. In particular, this fact asserts that the space of basic semialgebraic sets is invariant under the restricted preimage map
${F}^{-1}_{\mathbb{X}}$. A simple constructive argument demonstrating 
this point is as follows. Fix any arbitrary basic semialgebraic set $X$ in $\mathbb{R}^n$
\begin{equation}
\label{eq:03.04}
{X}:=\left\{ x\in\mathbb{R}^{n} \ : \ \varphi(x)\geq 0 \right\},
\end{equation}  
where $\varphi\in\mathbb{R}^{p}[x]$ is known exactly. The value of the
restricted preimage map ${F}_{\mathbb{X}}^{-1}(\cdot)$ at $X$ can be directly evaluated as follows
\begin{align} 
{F}_{\mathbb{X}}^{-1}({X}) &=  \left\{   x\in\mathbb{R}^{n}  \ : \ \varphi_{0}(x) \geq 0, \ \varphi( f(x) ) \geq 0  \right\} \nonumber\\
\label{eq:03.05}
&= \left\{   x\in\mathbb{R}^{n}  \ : \ \varphi^{+}( x ) \geq 0 \right\},
\end{align}
where $\varphi^{+}$ is given by
\begin{equation}
\label{eq:03.06}
\varphi^{+}(x) = 
\begin{pmatrix}
\varphi_{0}(x) \\
\varphi( f(x) )
\end{pmatrix}.
\end{equation}
By Assumption~\ref{ass:02.01} and the above postulated hypothesis on ${X}$, $f$ and $\varphi$ are polynomials. Thus, $\varphi \circ f$ is also a polynomial. Since $\varphi_{0}$ is also a polynomial by Assumption~\ref{ass:02.02}, it follows that $\varphi^{+}$ is polynomial; In particular, $\varphi^{+}\in\mathbb{R}^{p+m}[x]$. Consequently, as asserted, the value ${F}^{-1}_{\mathbb{X}}({X})$ of the restricted preimage map ${F}^{-1}_{\mathbb{X}}$ evaluated at $X$  is itself a basic semialgebraic set. 
\begin{thm}
\label{thm:03.01} 
Suppose Assumptions~\ref{ass:02.01} and~\ref{ass:02.02} hold.  ${F}^{-1}_{\mathbb{X}}({X})$ 
is a basic semialgebraic set in $\mathbb{R}^n$ for any basic semialgebraic set $X$ in $\mathbb{R}^n$.
\end{thm}

A direct use of induction and Theorem~\ref{thm:03.01} yields the relevant structural properties of the maximal trajectory, specified by~\eqref{eq:03.03}, of the set--dynamics of~\eqref{eq:03.01}. The related properties are summarized in the following theorem.
\begin{thm}
\label{thm:03.02}
Suppose Assumptions~\ref{ass:02.01} and~\ref{ass:02.02} hold. The sequence 
$\{ {X}_{k}\}_{k\in\mathbb{N}}$ of sets in $\mathbb{R}^{n} $ generated by set iteration~\eqref{eq:03.03} is
\begin{enumerate}[label=(\alph*),ref=(\alph*)]
\item \label{thm:03.02.a} Monotonically nonincreasing, i.e. 
 $$\forall k\in\mathbb{N},\ {X}_{k+1}\subseteq {X}_{k}.$$
\end{enumerate}
Moreover, each term ${X}_{k}$ of this sequence is
\begin{enumerate}[label=(\alph*),ref=(\alph*),resume]
\item \label{thm:03.02.b}   Compact (possibly empty) and
\item \label{thm:03.02.c}       Basic semialgebraic.
\end{enumerate}
\end{thm}

\subsection{Nonemptiness}
\label{sec:03.02}

The nonemptiness of the limit  and terms of the maximal trajectory of the set--dynamics of~\eqref{eq:03.01} provides necessary and sufficient conditions for satisfaction of state constraints by original state dynamics over, respectively, the infinite and finite time horizons. In particular, a state trajectory generated by original state dynamics of~\eqref{eq:02.01} satisfies state constraints of~\eqref{eq:02.02} over the horizon $\{0,1,2,\ldots,k\}$ if and only if the $k^\text{th}$ term $X_k$ (of the set sequence 
$\{ {X}_{k}\}_{k\in\mathbb{N}}$ generated by set iteration~\eqref{eq:03.03}) is nonempty and the initial state $x_0$ belongs to it. Likewise, a state trajectory generated by original state dynamics of~\eqref{eq:02.01} satisfies state constraints of~\eqref{eq:02.02} indefinitely long (i.e. over the infinite horizon $\{0,1,2,\ldots,k,\ldots\}$) if and only if the limit of the maximal trajectory of the set--dynamics of~\eqref{eq:03.01} is nonempty and the initial state $x_0$ belong to it. The case in which there is an empty term $X_k$ and, thus, the limit is empty is important but trivial as it yields a conclusion that no trajectory exists that satisfies state constraint over the horizon $\{0,1,2,\ldots,k\}$ and, thus, over the infinite horizon $\{0,1,2,\ldots,k,\ldots\}$. The highly desired nonemptiness can be guaranteed under a rather mild condition, which is frequently encountered in the usual setting of polynomial systems.

\begin{assum}
\label{ass:03.01}
The function $f$ admits at least one fixed point $\bar{x}=f(\bar{x})$ in 
$\mathbb{X}$.
\end{assum}

The main ramification of this additional, and clearly natural, assumption is the following strengthening of Theorem~\ref{thm:03.02}.
\begin{cor}
\label{cor:03.01}
Suppose Assumptions~\ref{ass:02.01},~\ref{ass:02.02}, and~\ref{ass:03.01} 
hold. Then, in addition to the properties stated in 
Theorem~\ref{thm:03.02}, each set ${X}_{k}$ is nonempty.
\end{cor}

\subsection{The Maximal Positively Invariant Set}
\label{sec:03.03}

It is reasonably well known~\cite{bertsekas:1972,gilbert:tan:1991,aubin:1991,kolmanovsky:gilbert:1998,blanchini:miani:2008}  that the maximal positively invariant set $\Omega_{\infty}$ for dynamics of~\eqref{eq:02.01} and constraint set of~\eqref{eq:02.02} is the unique maximal fixed point of the map $F^{-1}_{\mathbb{X}}$. In fact, it is equal to the limit of the nonexpansive sequence $\{{X}_{k}\}_{k\in\mathbb{N}}$, 
i.e.
\begin{equation}
\label{eq:03.07}
\Omega_{\infty} = \bigcap\limits_{k\geq 0} {X}_{k}.
\end{equation} 
\emph{Since the class of basic semialgebraic sets is not closed under infinite intersections, the only topological properties that can be guaranteed for sure for $\Omega_{\infty}$, under our mild conditions, are its compactness and nonemptiness.}

\subsection{Sufficient Conditions for Finite Time Termination}
\label{sec:03.04}

A term ${X}_{k}$ of the sequence generated through set recursion~\eqref{eq:03.03} is a positively invariant set for 
dynamics of~\eqref{eq:02.01} and constraint set of~\eqref{eq:02.02} if and only if ${X}_{k}\subseteq{X}_{k+1}$ in which case, in fact, it holds that ${X}_{k}={X}_{k+1}$. Hence, since the set sequence $\{{X}_{k}\}_{k\in\mathbb{N}}$ generated by~\eqref{eq:03.03} is nonexpansive, a necessary and sufficient condition for any of its terms ${X}_{k}$ to be the maximal positively invariant set  for dynamics of~\eqref{eq:02.01} and constraint set of~\eqref{eq:02.02} reads as
\begin{equation}
\label{eq:03.08}
{X}_{k+1} = {F}^{-1}_{\mathbb{X}}({X}_k) = {X}_{k}.
\end{equation}
Without further assumptions on the dynamics and constraints, the finite time termination of the set  iteration~\eqref{eq:03.03} is not guaranteed.  Fortunately, sufficient conditions on the state transition map $f$ and the constraint set $\mathbb{X}$ guaranteeing finite time termination are relatively mild and practically applicable.

\begin{assum}
\label{ass:03.02}
The (not necessarily unique) fixed point $\bar{x}$ of the state transition map $f$ is locally asymptotically stable with the basin of attraction $\Psi$. Moreover, $\bar{x}$ is contained in the interior of the constraint set $\mathbb{X}$, and $\mathbb{X}\subseteq \Psi$.
\end{assum}

\begin{thm}
\label{thm:03.03}
Suppose Assumptions~\ref{ass:02.01},~\ref{ass:02.02}, and~\ref{ass:03.02} hold. Then there exists a finite $k\in\mathbb{N}$ such that the set ${X}_{k}$ generated through the set iteration~\eqref{eq:03.03} satisfies~\eqref{eq:03.08} and, thus, it is the maximal positively invariant set for dynamics of~\eqref{eq:02.01} and constraint
set of~\eqref{eq:02.02}. 
\end{thm}

Theorem~\ref{thm:03.03} identifies the computationally most relevant case.
\begin{cor}
\label{cor:03.02}
Suppose Assumptions~\ref{ass:02.01},~\ref{ass:02.02}, and~\ref{ass:03.02} 
hold. Then, the maximal positively invariant set $\Omega_\infty$ for dynamics of~\eqref{eq:02.01} and constraint
set of~\eqref{eq:02.02} is a (nonempty) compact and basic semialgebraic set.
\end{cor}

\section{Computational Aspects}
\label{sec:04}

The set recursion~\eqref{eq:03.03} together with the termination criterion~\eqref{eq:03.08} suggests a prototype algorithm for computing the maximal positively invariant set. From a computational perspective, there are clearly two steps. The first one is the construction of the explicit of implicit form of the map $\varphi_{k+1}$ defining the set ${X}_{k+1}$. The second one is the verification of the termination criterion ${X}_{k+1} = {X}_{k}$. 

\subsection{Explicit and Implicit Representations of Sets $X_k$}
\label{sec:04.01}

A direct iteration of~\eqref{eq:03.06} shows that, at every step of set recursion~\eqref{eq:03.03}, the set $X_k$ takes form
\begin{equation}
\label{eq:04.01}
{X}_{k} := \{ x\in\mathbb{R}^{n} \ : \ \varphi_{k}(x)\geq 0  \},
\end{equation}
with, in the worst case, $\varphi_{k}\in\mathbb{R}^{m(k+1)}[x]$, and where
\begin{equation}
\label{eq:04.02}
\varphi_{k} = 
\begin{pmatrix}
\partial\varphi_{0} \\
\vdots \\
\partial\varphi_{k-1}\\
\partial\varphi_{k}
\end{pmatrix}
\end{equation}
 and for each $k\in \mathbb{N}$, the partial update maps $\partial \varphi_k$ are given by
\begin{equation}
\label{eq:04.03}
\partial\varphi_k=\varphi_0 \circ f^{k}
\end{equation} 
with $f^0(x):=x$ and, in turn, $\partial \varphi_0=\varphi_0$. The explicit representation of the sets $X_k$ requires polynomial compositions~\eqref{eq:04.03} to be performed explicitly in order to obtain explicit forms of the maps $\varphi_k$. In view of~\eqref{eq:04.02} and~\eqref{eq:04.03}, the number of polynomials necessary to describe the sets $X_k$ increases linearly with $k$, while the degree of the associated polynomials increases exponentially. 

The implicit representation of the sets $X_k$ does not require polynomial compositions~\eqref{eq:04.03} to be constructed, and instead it employs implicit forms of the maps $\varphi_k$. More precisely, in this case, the set $X_{k}$ takes the form
\begin{equation}
\label{eq:04.04}
{X}_{k} := \left\{ x\in\mathbb{R}^{n} \ : \ 
\begin{aligned}
&\varphi_{0}(x)\geq 0,\ \varphi_{0}(f(x))\geq 0,\ldots, \\ 
&\varphi_{0}(f^{k-1}(x))\geq 0,\ \varphi_{0}(f^k(x))\geq 0  
\end{aligned}
\right\},
\end{equation}
where the polynomial compositions $\varphi_0\circ f^k$ are not performed directly, but are taken as the implicit forms of the partial update maps $\partial\varphi_k$. In this case, the number of implicit forms of the partial maps $\partial\varphi_k$ necessary to describe sets $X_k$ increases linearly with $k$, while there is no exponential increase of the degree of the related polynomials due to the use of the corresponding implicit forms. 

\subsection{Verification of Finite Time Termination}
\label{sec:04.02}

Since the iterates satisfy  ${X}_{k+1}\subseteq{X}_{k}$ for every $k\in\mathbb{N}$, the verification of the termination criterion reduces to verifying the reverse inclusion  ${X}_{k}\subseteq{X}_{k+1}$. Now, in light of the structure of ${X}_{k+1}$ specified in~\eqref{eq:04.02}, the set inclusion ${X}_{k}\subseteq {X}_{k+1}$ holds true if and only if each polynomial $\partial\varphi_{(k+1,i)}$ with $i\in\{1,\ldots,m\}$ is nonnegative over ${X}_{k}$. This condition can be checked by solving $m$ polynomial optimization problems specified, for all $i\in \{1,\ldots,m\}$, by
\begin{equation}
\label{eq:04.05}
d_{i}  := \min_{x} \{\partial\varphi_{(k+1,i)}(x)\ :\ x\in X_k\},
\end{equation}
as the termination condition~\eqref{eq:03.08} holds true if and only if 
\begin{equation}
\label{eq:04.06}
\delta := \min_{i}\{d_{i} \ : \ i\in\{1,\ldots,m\} \}\geq 0.
\end{equation}

\subsection{Prototype Algorithm}
\label{sec:04.03}

When working with implicit representations of the sets $X_k$, the implicit representation of the maximal positively invariant set $\Omega_\infty$ is direct to construct as soon as an integer $k\in \mathbb{N}$ verifying finite time termination is detected. The latter chore is the main and nontrivial computational burden in this case.

\begin{algorithmn}{\label{alg:04.01} The prototype algorithm for the computation of the explicit form of the maximal positively invariant set.}
\vspace{0.4em}
{\footnotesize
\begin{algorithmic}[2]
\REQUIRE State transition map $f\in\mathbb{R}^{n}[x]$ $\bullet$ Constraint map
$\varphi_{0}\in\mathbb{R}^{m}[x]$ $\bullet$ Maximum number of iterations 
$k_\text{max}\geq 0$ $\bullet$ Tolerance $\varepsilon\geq 0$. 
\\[0.2em]
\STATE Set $k \leftarrow 0$,  $\delta \leftarrow -1$,  $\varphi_{k} \leftarrow \varphi_{0}$, and $\partial\varphi_{0} \leftarrow \varphi_{0}$
\vspace{0.2em}
\WHILE{$\delta < \varepsilon$ \AND $k<k_\text{max}$}
\vspace{0.4em}
\STATE{$
\partial\varphi_{k+1} \leftarrow 
\varphi_{0} \circ f^{k+1}
$}
\vspace{0.4em}
\FORALL{$i\in\{1,\ldots,m\}$}
\vspace{0.4em}
\STATE{
$
d_{i} \leftarrow \min_{x}\left\{ \partial\varphi_{(k+1,i)}(x) \ : \ 
\varphi_{k}(x) \geq 0 \right\}
$}
\vspace{0.4em}
\ENDFOR
\vspace{0.4em}
\STATE $\delta \leftarrow \min\{ d_{1},\ldots,d_{m} \} $ 
\vspace{0.4em}
\IF{$\delta < \varepsilon$} 
\vspace{0.4em}
\STATE{$k \leftarrow k+1$ \AND
$
\varphi_{k+1} \leftarrow 
\begin{pmatrix}
\varphi_{k} \\[0.2em]
\partial\varphi_{k+1}
\end{pmatrix}
$}
\vspace{0.4em}
\STATE Optional: Complexity reduction
\vspace{0.4em}
\ENDIF 
\vspace{0.4em}
\ENDWHILE
\vspace{0.4em}
\IF{$\delta\geq\varepsilon$}
\vspace{0.4em}
\RETURN Map $\varphi_{k}\in\mathbb{R}^{m(k+1)}[x]$ defining $X_k$ and $\Omega_\infty$
\vspace{0.4em}
\ELSE
\vspace{0.4em}
\RETURN UNSUCCESFUL
\vspace{0.4em}
\ENDIF
\end{algorithmic}
}
\vspace{0.4em}
\end{algorithmn}

When working with explicit representations of the sets $X_k$, the computation of the explicit representation of the maximal positively invariant set $\Omega_\infty$ requires execution of a number of polynomial compositions in addition to verification of finite time termination. Verification of finite time termination is computationally more convenient with explicit forms of the involved sets. But this potential computational convenience is offset by the requirements to perform the corresponding polynomial compositions. When implementing the above outlined procedures, it is also important to add another termination criterion. For example, one may want to repeat the outlined steps until a maximum number of iterations $k_\text{max}$ has been reached. When the dynamics are locally exponentially stable, bounds on the maximum number of steps the set recursion~\eqref{eq:03.03} may take prior to converging can be computed using relatively direct arguments based on Lyapunov exponents~\cite{rakovic:fiacchini:2008ifac:b}. (Naturally, such bounds can be also directly employed to construct implicit representation of the maximal positively invariant set.) Also, it might be worthwhile realizing that the case when $X_k$ is empty for some $k$ leads to a conclusion that $\Omega_{\infty}=X_{k+1}=X_{k}=\emptyset$. This case is omitted from the prototype algorithm since it is irrelevant from a practical point of view. A summary of the basic steps to compute the explicit form of the maximal positively invariant set for dynamics of~\eqref{eq:02.01} and constraint set of~\eqref{eq:02.02} is provided in Algorithm~\ref{alg:04.01}. 

\subsection{Miscellaneous Computational Remarks}
\label{sec:04.04}

Checking boundedness (and hence compactness) of $X_{0}$ can be done by solving 
\begin{equation}
\label{eq:04.07}
\max_{x}\left\{c_{j}^\intercal x \ : \ \varphi_{0}(x) \geq 0 \right\},
\end{equation}
where the vectors $c_{j}\in\mathbb{R}^{n}$ are the normal vectors of the $n$--dimensional unit simplex. Testing nonemptiness of the constraint  set $X_0$ can be done by solving
\begin{equation}
\label{eq:04.08}
\max_{x}\{ 0 \ : \ \varphi_{0}(x) \geq 0  \}.
\end{equation}
Detecting a constraint admissible fixed point can be done by solving a polynomial optimization of the following form
\begin{equation}
\label{eq:04.09}
\min_{x}\left\{\Vert f(x) - x\Vert^{2}_{2} \ : \ \varphi_{0}(x) \geq 0 \right\}.
\end{equation}
Verifying local asymptotic stability of the fixed point is more involved, but sum--of--squares certificates can be computed by making use of recent techniques developed in~\cite{ahmadi:2011}.

 An optional step in the above algorithm is removal of the number of redundant polynomial constraints. Suppose that a polynomial map $\varphi_{k+1}\in\mathbb{R}^{\ell}[x]$ is constructed at the end of the $k^\text{th}$ iteration of the prototype algorithm. This removal can be achieved by solving at each step $\ell$ polynomial optimization problems specified, for all $i\in\{1,\ldots,\ell\}$, by
\begin{equation}
\label{eq:04.10}
r_{i} = \min_{x} \left\{\varphi_{(k+1,i)}(x) \ : \ 
\begin{aligned}
&\forall j\in  \{1,\dots,\ell\}\setminus\{i\},\\
&\varphi_{(k+1,j)}(x) \geq 0
\end{aligned}
\right\}
\end{equation}
and removing the constraints $\varphi_{(k+1,i)}(x)$ for which $r_{i}\geq 0$. Although for the class of linear (and, in general, convex) polynomials the reduction step is clearly beneficial, for the general case the benefits may not be clear since the redundant constraints often help to tighten the relaxations used to solve nonconvex global optimization problems~\cite{shor:1998}. Finally, the underlying polynomial opimization is computationally hard, but it can be implemented by employing the existing third--party MATLAB\textsuperscript{\textregistered} toolboxes {\tt SOSTOOLS}~\cite{sostools} and {\tt GloptiPoly}~\cite{henrion:lassere:lofberg:2009}.  

\section{Numerical Example}
\label{sec:05}

\begin{figure}
\centering  
\includegraphics[width=0.86\columnwidth]{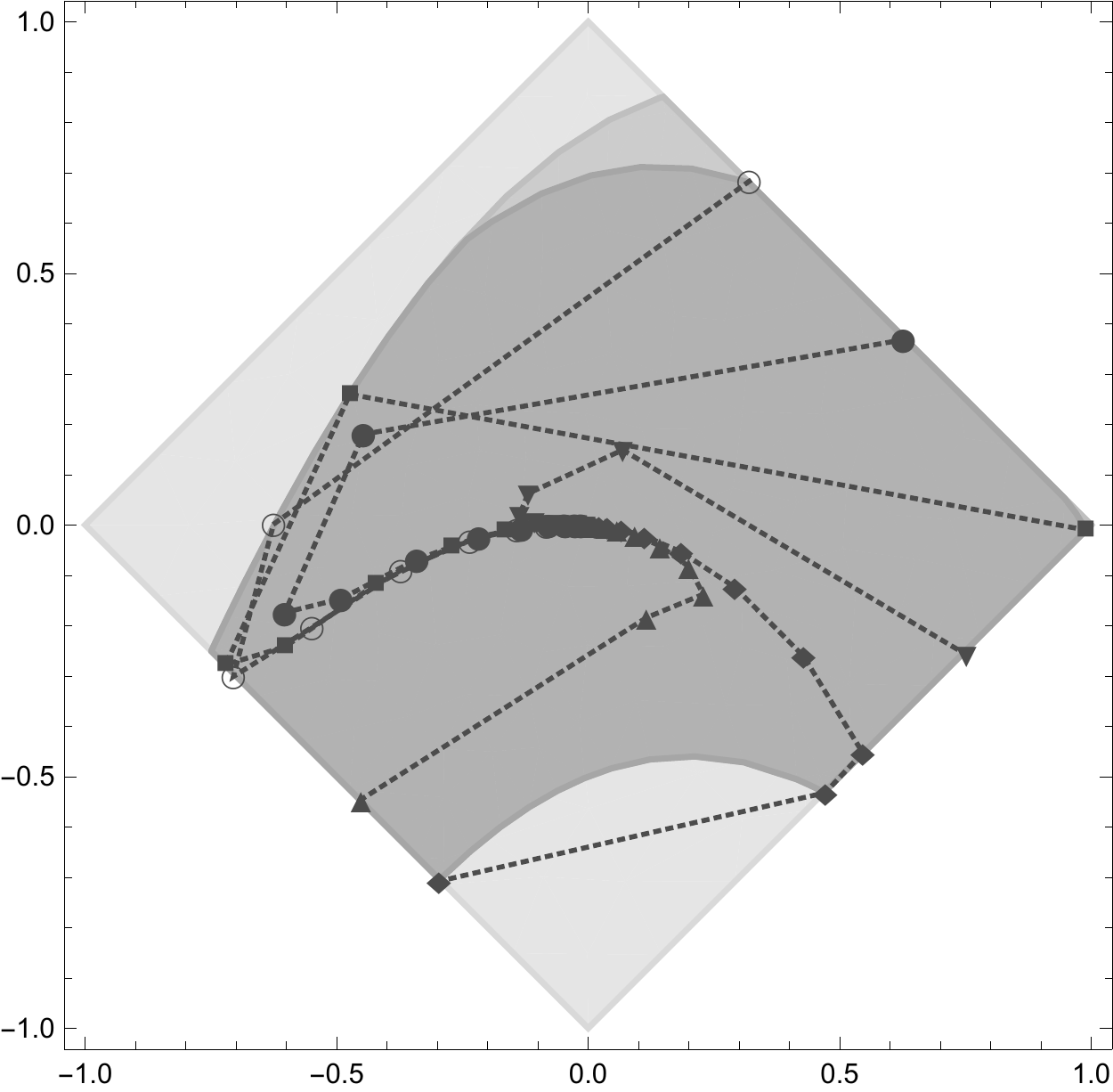}
\vspace{1em}

\includegraphics[width=0.86\columnwidth]{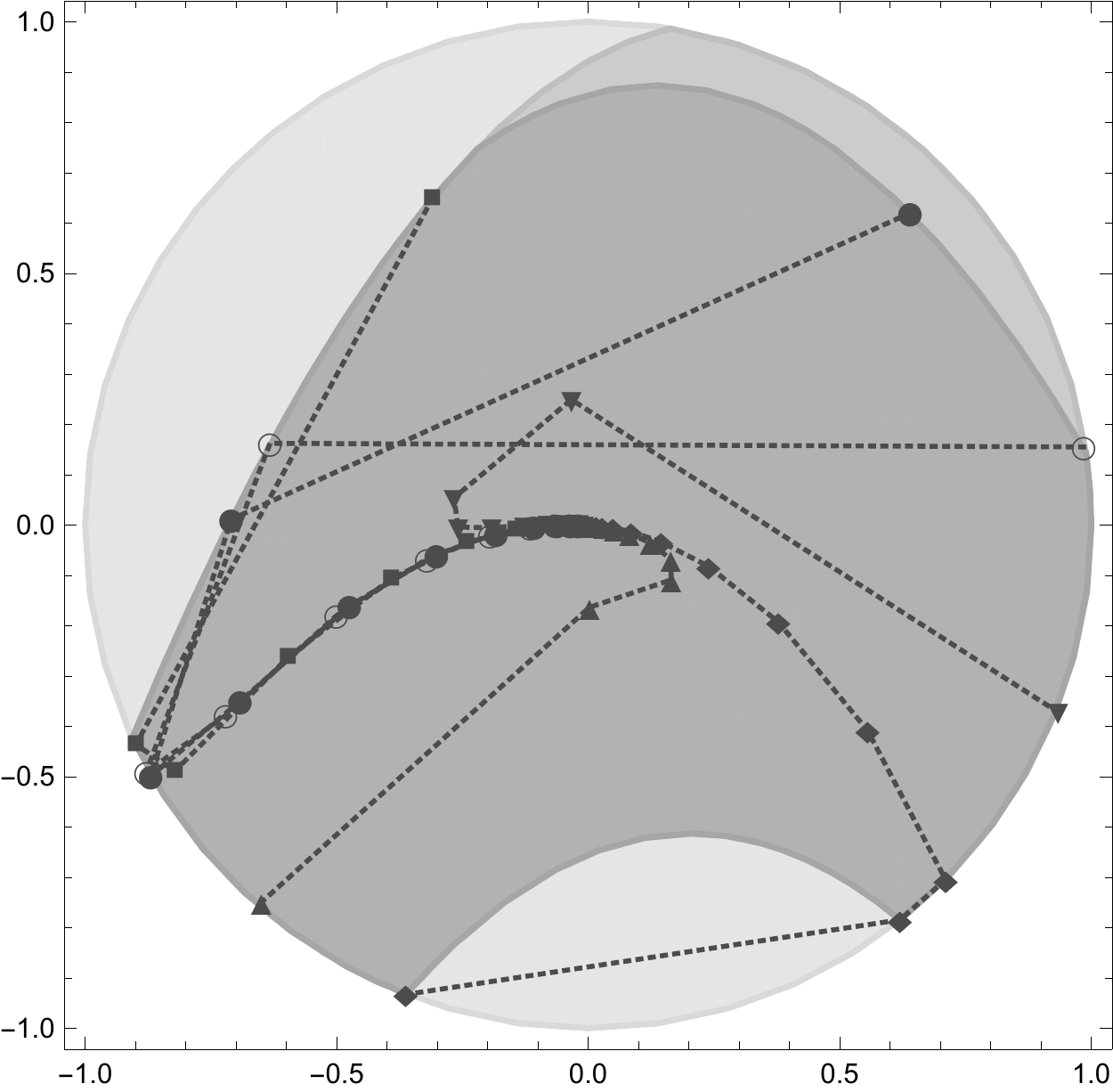} 
\vspace{1em}

\includegraphics[width=0.86\columnwidth]{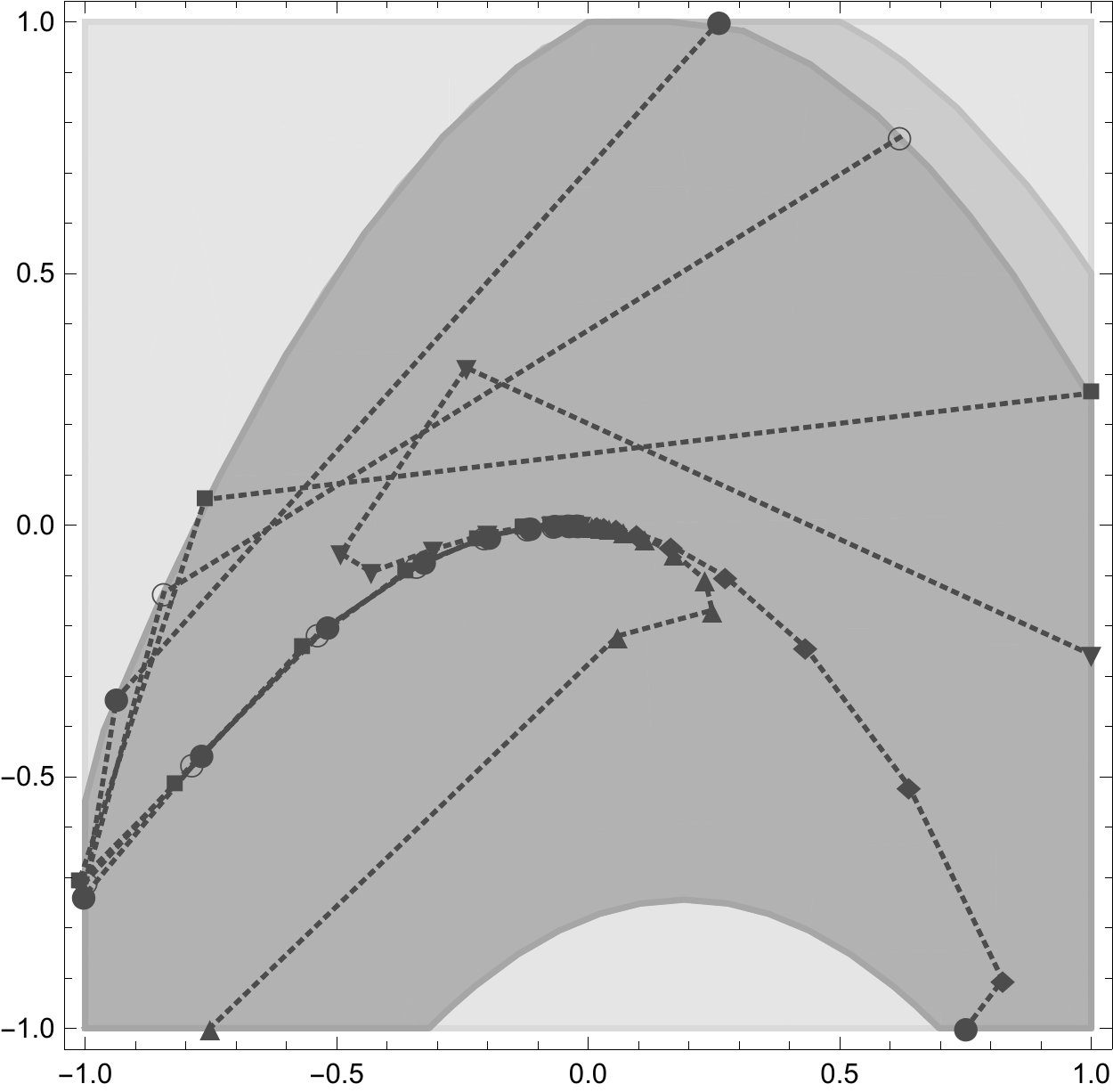}
\caption{\label{fig:05.01}The maximal positively invariant set of  dynamics of~\eqref{eq:05.01} constrained
by the 1-norm ball (top), 2-norm ball (middle), and infinity--norm ball (bottom).
Color key: ${\color{gray90}\blacksquare}\; X_{0} = \mathbb{X}$, 
${\color{gray80}\blacksquare}\; X_{1}$, and ${\color{gray70}\blacksquare}\; 
X_{2} = \Omega_{\infty}$.   } 
\end{figure}

Consider the following polynomial dynamics
\begin{equation}
\label{eq:05.01}
x^{+} = 
\begin{pmatrix}
\frac{1}{2}x_{1} - x_{1}^{2} - x_{2} \\
\frac{1}{2}x_{1}^{2} + \frac{1}{2}x_{2}-(\frac{1}{2}x_{1}-x_{1}^{2}-x_{2}^2) 
\end{pmatrix},
\end{equation}
where $x=(x_1,x_2)\in\mathbb{R}^{2}$ and $x^{+}=(x_1^+,x_2^+)\in\mathbb{R}^{2}$ are, respectively, the current and successor states. As shown in~\cite{ahmadi:2011}, the origin is a globally asymptotically stable equilibrium point for the system. In order to illustrate the effect of the constraint set on the maximal positively invariant set, we consider the cases where the constraint set is given by 1) the vector 1--norm ball
\begin{equation}
\label{eq:05.02}
\mathbb{X} = \left\{ \ x\in\mathbb{R}^{2} \ : \
\begin{alignedat}{2}
&1 - (x_{1} + x_{2}) &&\geq 0 \\
&1 - (-x_{1} + x_{2}) &&\geq 0 \\
&1 - (-x_{1} - x_{2}) &&\geq 0 \\ 
&1 - (x_{1} - x_{2}) &&\geq 0
\end{alignedat} \ \right\},
\end{equation}
2) the Euclidean vector norm ball 
\begin{equation}
\label{eq:05.03}
\mathbb{X} = \{ x\in\mathbb{R}^{2} \ : \ 1-(x_{1}^{2} + x_{2}^{2} ) \geq 0 \},
\end{equation}
and 3) the vector $\infty$--norm ball 
\begin{equation}
\label{eq:05.04}
\mathbb{X} = \left\{ \ x\in\mathbb{R}^{2} \ : \
\begin{alignedat}{2}
&1 - x_{1} &&\geq 0 \\
&1 - x_{2} &&\geq 0 \\
&1 + x_{1} &&\geq 0 \\ 
&1 + x_{2} &&\geq 0
\end{alignedat} \ \right\}.
\end{equation}
In all instances, the algorithm 
terminated in the third iteration, providing a vector of polynomials describing
the maximal positively invariant set. The iterates $\mathbb{X}={X}_{0}$, ${X}_{1}$, and, ${X}_{2}=\Omega_{\infty}$
are shown in Figure~\ref{fig:05.01} together with sample trajectories depicted by points (different markers different trajectories) and their linear interpolations depicted by dotted lines. The figure was generated in Wolfram Mathematica $10$.

\section{Conclusion}
\label{sec:06}

The structural properties of the maximal positively invariant set for polynomial discrete time dynamics subject to basic semialgebraic constraints were considered.   The basic semialgebraic structure of the maximal positively invariant set was established under relatively mild assumptions. A prototype algorithm for the computation of the maximal positively invariant set, based on the existing tools for polynomial optimization, was  provided and illustrated by an example.

\begin{ack}                               
The authors would like to thank the Texas A\&M Energy Institute, 
and in particular Prof. Efstratios N. Pistikopoulos for providing 
financial support. 
\end{ack}

\bibliographystyle{plain}        
\bibliography{mpis-ps-rv}           

\begin{thebibliography}{10}

\bibitem{ahmadi:2011}
A.~A. Ahmadi.
\newblock {\em {Algebraic Relaxations and Hardness Results in Polynomial
  Optimization and Lyapunov Analysis}}.
\newblock PhD thesis, Massachisetts Institute of Technology, 2011.

\bibitem{artstein:rakovic:2008}
Z.~Artstein and S.~V. Rakovi\'{c}.
\newblock {Feedback and Invariance under Uncertainty via Set Iterates}.
\newblock {\em Automatica}, 44(2):520--525, 2008.

\bibitem{artstein:rakovic:2011}
Z.~Artstein and S.~V. Rakovi\'{c}.
\newblock {Set Invariance Under Output Feedback : A Set--Dynamics Approach}.
\newblock {\em {International Journal of Systems Science}}, 42(4):539--555,
  2011.

\bibitem{aubin:1991}
J.~P. Aubin.
\newblock {\em {Viability Theory}}.
\newblock Systems \& {C}ontrol: {F}oundations \& {A}pplications. Birkhauser,
  Boston, {B}asel, {B}erlin, 1991.

\bibitem{bertsekas:1972}
D.~P. Bertsekas.
\newblock {Infinite--time Reachability of State--space Regions by Using
  Feedback Control}.
\newblock {\em IEEE Trans. Automatic Control}, 17(5):604--613, 1972.

\bibitem{blanchini:miani:2008}
F.~Blanchini and S.~Miani.
\newblock {\em {Set--Theoretic Methods in Control}}.
\newblock Systems \& {C}ontrol: {F}oundations \& {A}pplications.
  Birkh{\"a}user, Boston, {B}asel, {B}erlin, 2008.

\bibitem{gilbert:tan:1991}
E.~G. Gilbert and K.~T. Tan.
\newblock {Linear Systems with State and Control Constraints: the Theory and
  Application of Maximal Output Admissible Sets}.
\newblock {\em IEEE Transactions on Automatic Control}, 36(9):1008--1020, 1991.

\bibitem{henrion:lassere:lofberg:2009}
D.~Henrion, J.~B. Lasserre, and J.~L{\"o}fberg.
\newblock {GloptiPoly 3: Moments, Optimization and Semidefinite Programming}.
\newblock {\em Optimization Methods and Software}, 24(4--5):761--779, 2009.

\bibitem{kolmanovsky:gilbert:1998}
I.~V. Kolmanovsky and E.~G. Gilbert.
\newblock {Theory and Computation of Disturbance Invariant Sets for Discrete
  Time Linear Systems}.
\newblock {\em Mathematical Problems in Engineering: Theory, Methods and
  Applications}, 4:317--367, 1998.

\bibitem{korda:henrion:jones:2014}
M.~Korda, D.~Henrion, and C.~Jones.
\newblock {Convex Computation of the Maximum Controlled Invariant Set For
  Polynomial Control Systems}.
\newblock {\em {SIAM Journal on Control and Optimization}}, 52(5):2944--2969,
  2014.

\bibitem{li:liu:2016cacsd}
Y.~Li and J.~Liu.
\newblock {Computing Maximal Invariant Sets for switched Nonlinear Systems}.
\newblock In {\em Proceedings of the 2016 IEEE Conference on Computer Aided
  Control System Design}, Buenos Aires, Argentina, 2016.

\bibitem{liu:zhan:zhao:2011emsoft}
J.~Liu, N.~Zhan, and H.~Zhao.
\newblock {Computing Semi--Algebraic Invariants for Polynomial Dynamical
  Systems}.
\newblock In {\em Proceedings of the ninth ACM international conference on
  Embedded software EMSOFT 2011}, New York, NY, USA, 2011.

\bibitem{olaru:dedona:seron:stoican:2011}
S.~Olaru, J.~A. {De Don\'{a}}, M.~M. Seron, and F.~Stoican.
\newblock {Positive Invariant Sets for Fault Tolerant Multisensor Control
  Schemes}.
\newblock {\em {International Journal of Control}}, 83(12):2622--2640, 2010.

\bibitem{sostools}
A.~Papachristodoulou, J.~Anderson, G.~Valmorbida, S.~Prajna, P.~Seiler, and
  P.~A. Parrilo.
\newblock {\em {SOSTOOLS: Sum of Squares Optimization Toolbox for MATLAB}}.
\newblock \texttt{http://arxiv.org/abs/1310.4716}, 2013.
\newblock Available from \texttt{http://www.eng.ox.ac.uk/control/sostools},
  \texttt{http://www.cds.caltech.edu/sostools} and
  \texttt{http://www.mit.edu/\~{}parrilo/sostools}.

\bibitem{parrilo:2003}
P.~Parrilo.
\newblock {Semidefinite Programming Relaxations for Semialgebraic Problems}.
\newblock {\em Mathematical Programming}, 96(2):293--320, 2003.

\bibitem{rakovic:fiacchini:2008ifac:b}
S.~V. Rakovi\'{c} and M.~Fiacchini.
\newblock {Invariant Approximations of the Maximal Invariant Set of
  ``Encircling the Square''}.
\newblock In {\em Proceedings of the 17th IFAC World Congress IFAC 2008},
  Seoul, Korea, 2008.

\bibitem{rawlings:mayne:2009}
J.~B. Rawlings and D.~Q. Mayne.
\newblock {\em {Model Predictive Control: Theory and Design}}.
\newblock {N}ob {H}ill {P}ublishing, Madison, 2009.

\bibitem{rudin:1964}
W.~Rudin.
\newblock {\em {Principles of Mathematical Analysis}}.
\newblock McGraw-Hill New York, 1953.

\bibitem{shor:1998}
N.~Z. Shor.
\newblock {Role of Redundant Constraints for Improving Dual Bounds in
  Polynomial Optimization Problems}.
\newblock {\em Cybernetics and Systems Analysis C/C OF Kibernetika i Sistemnyi
  Analiz}, 34:564--576, 1998.

\end{thebibliography}
\appendix
\section*{Appendix A: Proof of Theorem~\ref{thm:03.02}}
\paragraph*{\ref{thm:03.02.a}} Suppose that ${X}_{k+1}\subseteq{X}_{k}$ for some $k\in \mathbb{N}$. Then, 
$X_{k+2}={F}^{-1}_{\mathbb{X}}( {X}_{k+1} )  \subseteq{F}^{-1}_{\mathbb{X}}( {X}_{k} )=X_{k+1}.$
Since $X_0\subseteq {F}^{-1}_{\mathbb{X}}( {X}_{0} )=X_1$, the claim follows by induction.

\paragraph*{\ref{thm:03.02.b}} The set $\{x\in\mathbb{R}^n\ :\ f(x)\in {X}_{k}\}$
is closed (possibly empty) so that the set ${F}^{-1}_{\mathbb{X}}( {X}_{k} )=\{x\in\mathbb{R}^n\ :\ f(x)\in {X}_{k}\}\bigcap \mathbb{X}$ is compact (possibly empty) because $\mathbb{X}$ is compact.  Since ${X}_{0}=\mathbb{X}$ is compact, the  claim also follows by induction. 

\paragraph*{\ref{thm:03.02.c}} First, ${X}_{0}=\mathbb{X}$ is a basic semialgebraic set. Second, by Theorem~\ref{thm:03.01}, ${X}_{k+1}={F}^{-1}_{\mathbb{X}}( {X}_{k} ) $ is basic semialgebraic whenever ${X}_{k}$ is so. The claim also follows by induction.

\section*{Appendix B: Proof of Corollary~\ref{cor:03.01}}
By Assumption~\ref{ass:03.01}, $\bar{x}=f(\bar{x})\in \mathbb{X}$. Thus, in fact, for all $k\in\mathbb{N}$, $\bar{x}=f^k(\bar{x})\in \mathbb{X}$. In turn, for all $k\in\mathbb{N}$, $\bar{x}\in{X}_{k}\not=\emptyset$.
\end{document}